\documentclass[10pt,reqno,a4paper,twoside,dvips]{article}

\usepackage{graphicx}
\usepackage{amsfonts}
\usepackage{amssymb}
\usepackage{graphicx}
\usepackage{amsthm}
\usepackage{amsmath}
\usepackage{eucal}
\usepackage{amsmath}
\newtheorem{teo}{Theorem}[section]
\newtheorem{defin}{Definition}
\newtheorem{lema}[teo]{Lemma}
\newtheorem{cor}{Corollary}[section]
\newtheorem{rema}{Remark}[section]
\newtheorem{exem}{Example}[section]
\newtheorem{prop}{Proposition}[section]
\newtheorem{pro}{P}

\newcommand{\ble}{\begin{lema}$\!\!\!\textrm{\bf{. \ \ \ }}$}
\newcommand{\ele}{\end{lema}}
\newcommand{\bde}{\begin{defin}$\!\!\!\textrm{\bf{.}}$}
\newcommand{\ede}{\end{defin}}
\newcommand{\bte}{\begin{teo}$\!\!\!\textrm{\bf{. }}$}
\newcommand{\ete}{\end{teo}}
\newcommand{\bob}{\begin{obs}$\!\!\!\textrm{\bf{.}}$}
\newcommand{\eob}{\end{obs}}
\newcommand{\bco}{\begin{cor}$\!\!\!\textrm{\bf{.}}$}
\newcommand{\eco}{\end{cor}}
\newcommand{\bcon}{\begin{con}$\!\!\!\textrm{\bf{.}}$}
\newcommand{\econ}{\end{con}}
\newcommand{\bre}{\begin{rema}$\!\!\!\textrm{\bf{. }}$}
\newcommand{\ere}{\end{rema}}
\newcommand{\bpr}{\begin{prop}$\!\!\!\textrm{\bf{.}}$}
\newcommand{\epr}{\end{prop}}
\newcommand{\bex}{\begin{exem}$\!\!\!\textrm{\bf{.}}$}
\newcommand{\eex}{\end{exem}}

\newcommand{\RR}{\mathbb R}

\newcommand{\bpp}{\begin{pro}$\!\!\!\textrm{\bf{.}}$}
\newcommand{\epp}{\end{pro}}

\newcommand{\be}{\begin{equation}}
\newcommand{\ee}{\end{equation}}
\newcommand{\ba}{\begin{array}}
\newcommand{\ea}{\end{array}}
\newcommand{\bi}{\begin{itemize}}
\newcommand{\ei}{\end{itemize}}
\newcommand{\bc}{\begin{center}}
\newcommand{\ec}{\end{center}}
\newcommand{\ov}{\overline}

\newcommand{\ds}{\displaystyle}

\numberwithin{equation}{section}

\begin{document}\thispagestyle{plain}

\title{Several applications of Cartwright-Field's inequality}
\author{Nicu\c{s}or
Minculete$^1$\footnote{E-mail: minculeten@yahoo.com } and   Shigeru
Furuichi$^2$\footnote{E-mail: furuichi@chs.nihon-u.ac.jp} }
\date{}
\maketitle \vspace{-1cm}
\begin{center}
$^1${\small \lq\lq Dimitrie Cantemir\rq\rq  University, Bra\c{s}ov,
500068, Rom{a}nia}\\
$^2${\small Department of Computer Science and System Analysis,}\\
{\small College of Humanities and Sciences, Nihon University,}\\
{\small 3-25-40, Sakurajyousui, Setagaya-ku, Tokyo, 156-8550, Japan}
\end{center}\vspace{1cm}

{\bf Abstract.} In this paper  we present several applications of
Cartwright-Field's inequality. Among these we found Young's
inequality, Bernoulli's inequality, the inequality between the
weighted power means, H\"{o}lder's inequality and Cauchy's
inequality. We give also two applications related to arithmetic
functions and to operator inequalities.

\vspace{3mm}

{\bf Keywords : } Cartwright-Field's inequality, Young's inequality,
Bernoulli's inequality, H\"{o}lder's inequality, arithmetic
function, operator inequality

\vspace{3mm} {\bf 2010 Mathematics Subject Classification : } 15A39,
15A45 and 26D15 \vspace{3mm}
\section{INTRODUCTION}
An important  result related to the improvement of the inequality
between arithmetic and geometric means (AM-GM) was obtained by D. I.
Cartwright and M. J. Field in \cite{CF}, which is given in the following
way: if $0<m=\min\{x_1,...,x_n\}$ and $M=\max\{x_1,...,x_n\}$, then
\begin{eqnarray}\label{ineq_1_1}
&\ds\frac{1}{2M}\sum_{i=1}^n\alpha_i\left(x_i-\sum_{k=1}^n \alpha_k
x_k\right)^2\leq \sum_{i=1}^n \alpha_i x_i-\prod_{i=1}^n
x_i^{\alpha_i}\leq&\notag\\
&\leq\ds \frac{1}{2m}\sum_{i=1}^n \alpha_i\left(x_i-\sum_{k=1}^n
\alpha_k x_k\right)^2,&
\end{eqnarray}  
where $\alpha_i>0$ for all $i=1...n$ and $\alpha_1+...+\alpha_n=1$.
For $n=2$, this inequality may be written as follows:
\begin{equation}\label{ineq_1_2}
\frac{\lambda(1-\lambda)}{2M}(a-b)^2\leq \lambda
a+(1-\lambda)b-a^\lambda b^{1-\lambda}\leq
\frac{\lambda(1-\lambda)}{2m}(a-b)^2,
\end{equation}
where $a,b>0, \ m=\min\{a,b\}, \ M=\max\{a,b\}$ and
$\lambda\in[0,1]$.
Since $\ds\frac{\lambda(1-\lambda)}{2M}(a-b)^2\geq0$, we deduce
Young's inequality (see \cite{HLP,NP})
\begin{equation}\label{ineq_1_3}
a^\lambda b^{1-\lambda} \leq \lambda a+(1-\lambda)b 
\end{equation}
Therefore, inequality (\ref{ineq_1_2}) is an improvement of Young's inequality
and at the same time gives a reverse inequality for the inequality
of Young.

In \cite{FM}, we presented two inequalities which give two different reverse inequalities for the Young's inequality, namely:
\begin{equation}\label{ineq_1_4}
0 \leq \lambda a+(1-\lambda) b -a^{\lambda}b^{1-\lambda} \leq
a^{\lambda}b^{1-\lambda}\exp\left\{\ds\frac{\lambda(1-\lambda)(a-b)^2}{m^2}\right\} -a^{\lambda}b^{1-\lambda} 
\end{equation}
and
\begin{equation}\label{ineq_1_5}
0 \leq \lambda a+(1-\lambda) b -a^{\lambda}b^{1-\lambda}  \leq
\lambda (1-\lambda)\left\{\log\left(\frac{a}{b}\right)\right\}^2M,
\end{equation}
where $a,b>0$, $m\equiv\min\{a,b\}$ and $M=\max\{a,b\}$  and
$\lambda\in[0,1]$.

\begin{rema}\label{rema_1_1}
The first inequality of (\ref{ineq_1_2}) clearly gives an improvement of the first inequality in (\ref{ineq_1_4}) and  (\ref{ineq_1_5}). 
For $0<a,b< 1$, we find the right hand side of the second inequality of (\ref{ineq_1_2}) gives tighter upper bound than that of (\ref{ineq_1_5}),
from the inequality $\frac{x-y}{\log x -\log y}<\frac{x+y}{2}$, for $x,y>0$.
For $a,b>1$, we find the right hand side of the second inequality of (\ref{ineq_1_5}) gives tighter upper bound than that of (\ref{ineq_1_2}),
from the inequality $\sqrt{xy}<\frac{x-y}{\log x -\log y}$, for $x,y>0$. 
In addition, we find the right hand side of the second inequality of  (\ref{ineq_1_2}) gives tighter upper bound than that of (\ref{ineq_1_4}) for $a,b>0$,
from $e^x > 1+ x$.
\end{rema}

Remark \ref{rema_1_1} supports the importance to study the inequality (\ref{ineq_1_2}) for several applications which will be given in the following sections.

\section{MAIN APPLICATIONS}

\begin{lema}
For $x>-1$ and $\lambda\in[0,1]$ there is the following inequality
\begin{equation}\label{ineq_2_1}
\frac{\lambda(1-\lambda)}{2M}x^2\leq \lambda x+1-(x+1)^\lambda \leq
\frac{\lambda(1-\lambda)}{2m}x^2,
\end{equation}
where $m=\min\{x+1,1\}$ and $M=\max\{x+1,1\}$.
\end{lema}

\noindent{\it Proof.} By replacing $\ds\frac{a}{b}$ of $t$ in
inequality (\ref{ineq_1_2}) we obtain the inequality
\begin{equation}\label{ineq_2_2}
\frac{\lambda(1-\lambda)}{2M}(t-1)^2\leq \lambda
t+1-\lambda-t^\lambda\leq \frac{\lambda(1-\lambda)}{2m}(t-1)^2,
\end{equation}
for all $t>0$ and $\lambda\in[0,1]$, where $m=\min\{t,1\}$ and
$M=\max\{t,1\}$.
Substituting $t=x+1$ in inequality (\ref{ineq_2_2}), we find the inequality
desired.

\hfill \qed

\begin{rema}
Inequality (\ref{ineq_2_1}) refines the inequality of Bernoulli, namely, for
$x>-1$ and $\lambda \in[0,1]$, we have
\begin{equation}
\lambda x+1\geq (x+1)^\lambda
\end{equation}
because $\ds\frac{\lambda(1-\lambda)}{2M}x^2\geq0$ in inequality
(\ref{ineq_2_1}).
\end{rema}
Next we will establish a refinement of the inequality between the
weighted power means, based on inequality (\ref{ineq_2_1}).

\begin{teo}
If $a_i>0, \ p_i>0, \ i=1...n, \ 0<r\leq s$,
$M_r(a,p)=\ds\left(\frac{\ds\sum_{i=1}^n p_i a_i^r}{\ds\sum_{i=1}^n
p_i}\right)^{1/r}$ and $M_s(a,p)=\left(\frac{\ds\sum_{i=1}^n p_i
a_i^s}{\ds\sum_{i=1}^n p_i}\right)^{1/s},$ then there is the
inequality
\begin{equation}\label{ineq_2_4}
\frac{A}{M}\leq[M_s(a,p)]^r-[M_r(a,p)]^r\leq\frac{A}{m},
\end{equation}
where
$$A=\frac{r(s-r)}{2s^2}[M_s(a,p)]^r\cdot\frac{\ds\sum_{i=1}^n p_i\left(\frac{a_i^s}{[M_s(a,p)]^s}-1\right)^2}{\ds\sum_{i=1}^n p_i},$$
$$m=\min_{i=\ov{1,n}}\left\{\frac{a_i^s}{[M_s(a,p)]^s},1\right\} \ {\rm and } \ M=\max_{i=\ov{1,n}}\left\{\frac{a_i^s}{[M_s(a,p)]^s},1\right\}.$$
\end{teo}

\noindent{\it Proof.} If $r=s$, then we have the equality in
relation
(\ref{ineq_2_4}). Let $r<s$.
 In inequality (\ref{ineq_2_2}) we consider
$t=\ds\frac{a_i^s}{[M_s(a,p)]^s}$ and $\lambda=\frac{r}{s}<1$, thus,
we deduce the inequality
\begin{eqnarray}
\ds\frac{r(s-r)}{2s^2M}\left(\frac{a_i^s}{[M_s(a,p)]^s}-1\right)^2 &\leq& \frac{r}{s}\frac{a_i^s}{[M_s(a,p)]^s}+1-\frac{r}{s}-\frac{a_i^r}{[M_s(a,p)]^r} \nonumber \\
&\leq&\ds\frac{r(s-r)}{2s^2m}\left(\ds\frac{a_i^s}{[M_s(a,p)]^s}-1\right)^2.\label{ineq_2_5}
\end{eqnarray}
Multiplying by $p_i$ in inequality (\ref{ineq_2_5}) and taking the sum for
$i=1...n$, we obtain the following inequality
\begin{eqnarray*}
\frac{r(s-r)}{2s^2M}\frac{\ds\sum_{i=1}^n p_i\left(\frac{a_i^s}{[M_s(a,p)]^s}-1\right)^2}{\ds\sum_{i=1}^n p_i} &\leq & 1-\left[\frac{M_r(a,p)}{M_s(a,p)}\right]^r \\
&\leq& \frac{r(s-r)}{2s^2 m}\frac{\ds\sum_{i=1}^n p_i\left(\ds\frac{a_i^s}{[M_s(a,p)]^s}-1\right)^2}{\ds\sum_{i=1}^n p_i},
\end{eqnarray*}
which is equivalent to the inequality of the statement.

\hfill \qed

\begin{rema}
Since $\ds\frac{A}{M}\geq0$ in inequality (\ref{ineq_2_4}), we find the
inequality between the weighted power means  \cite{HLP,NP},
\begin{equation}\label{ineq_2_6}
M_r(a,p)\leq M_s(a,p),
\end{equation}
for $0<r\leq s$.
The two means are equal if and only if $a_1=a_2=...=a_ n$.
\end{rema}

\begin{teo}
Let $p,q>1$ be real numbers satisfying
$\ds\frac{1}{p}+\frac{1}{q}=1$. If $a_i,b_i>0$ for all $i=1...n$,
then there is the following inequality
\begin{equation}\label{ineq_2_7}
\frac{A}{M}\leq\left(\sum_{i=1}^n
a_i^p\right)^{1/p}\left(\sum_{i=1}^n b_i^q\right)^{1/q}-\sum_{i=1}^n
a_ib_i\leq\frac{A}{m},
\end{equation}
where
$$A=\frac{1}{2pq}\left(\sum_{i=1}^n a_i^p\right)^{1/p}\left(\sum_{i=1}^n b_i^q\right)^{1/q}\sum_{i=1}^n \left(\frac{a_i^p}{\ds\sum_{i=1}^n a_i^p}-\frac{b_i^q}{\ds\sum_{i=1}^n b_i^q}\right)^2,$$
$$m=\min_{i=\ov{1,n}}\left\{\frac{a_i^p}{\ds\sum_{i=1}^n a_i^p}, \ \frac{b_i^q}{\ds\sum_{i=1}^n b_i^q}\right\} \ {\rm and } \ M=\max_{i=\ov{1,n}}\left\{\frac{a_i^p}{\ds\sum_{i=1}^n a_i^p}, \ \frac{b_i^q}{\ds\sum_{i=1}^n b_i^q}\right\}.$$
\end{teo}

\noindent{\it Proof.} By replacing $\lambda=\ds\frac{1}{p}, \
1-\lambda=\frac{1}{q}, \ a=\frac{a_i^p}{\ds\sum_{i=1}^n a_i^p}$;
$b=\ds\frac{b_i^q}{\ds\sum_{i=1}^n b_i^q}$ in inequality (\ref{ineq_1_2}) we
obtain the relation
$$\frac{1}{2pqM}\left(\frac{a_i^p}{\ds\sum_{i=1}^n a_i^p}-\frac{b_i^q}{\ds\sum_{i=1}^n b_i^q}\right)^2\leq \frac{a_i^p}{\ds p\sum_{i=1}^n a_i^p}+\frac{b_i^q}{q\ds\sum_{i=1}^m b_i^q}-\frac{a_i b_i}{\left(\ds\sum_{i=1}^na_i^p\right)^{1/p}\left(\ds\sum_{i=1}^n b_i^q\right)^{1/q}}$$
$$\leq \frac{1}{2pqm}\left(\frac{a_i^p}{\ds\sum_{i=1}^n a_i^p}-\frac{b_i^q}{\ds\sum_{i=1}^n b_i^q}\right)^2.$$
We observe that taking the sum for $i=1...n$ we deduce the
inequality of the statement.

\hfill \qed

\begin{rema}
\begin{itemize}
\item[(a)] H\"{o}lder's inequality is widely used in the theory of
inequalities and has the form \cite{HLP,NP}:
\begin{equation}\label{ineq_2_8}
\left(\sum_{i=1}^n a_i^p\right)^{1/p} \left(\ds\sum_{i=1}^n
b_i^q\right)^{1/q}\geq \sum_{i=1}^n a_i b_i.
\end{equation}
Because $\ds\frac{A}{M}\geq0$ in inequality (\ref{ineq_2_7}), we obtain a
proof of H\"{o}lder's inequality. It is easy to see that inequality
(\ref{ineq_2_7}) is a refinement of H\"{o}lder's inequality and contains a
reverse inequality for the inequality of H\"{o}lder.
\item[(b)] For $p=q=2$ in inequality (\ref{ineq_2_7}), we  have an improvement of
Cauchy's inequality
\begin{equation}\label{ineq_2_9}
\ds\sum_{i=1}^n a_i^2\sum_{i=1}^n b_i^2\geq \left(\sum_{i=1}^n a_i
b_i\right)^2,
\end{equation}
given by the inequality
\begin{equation}\label{ineq_2_10}
\frac{A^2}{M^2}+\frac{2A}{M}\sum_{i=1}^n a_ib_i\leq
\left(\sum_{i=1}^n a_i^2\right)\left(\sum_{i=1}^n
b_i^2\right)-\left(\sum_{i=1}^n a_ib_i\right)^2\leq
\frac{A^2}{m^2}+\frac{2A}{m}\sum_{i=1}^n a_i b_i,
\end{equation}
where
$$A=\frac{1}{8}\sqrt{\ds\left(\sum_{i=1}^n a_i^2\right)\left(\sum_{i=1}^n b_i^2\right)}\cdot \sum_{i=1}^n \left(\frac{a_i^2}{\ds\sum_{i=1}^n a_i^2}-\frac{b_i^2}{\ds\sum_{i=1}^n b_i^2}\right)^2,$$
$$m=\min_{i=\ov{1,n}}\left\{\frac{a_i^2}{\ds\sum_{i=1}^n a_i^2}, \ \frac{b_i^2}{\ds\sum_{i=1}^n b_i^2}\right\} \ {\rm and} \ M=\max_{i=\ov{1,n}}\left\{\frac{a_i^2}{\ds\sum_{i=1}^n a_i^2}, \ \frac{b_i^2}{\ds\sum_{i=1}^n b_i^2}\right\}.$$
The equality holds for $\ds\frac{a_1}{b_1}=...=\frac{a_n}{b_n}$.
\item[(c)] In \cite{Pop}, O. T. Pop gave Bergstr\"{o}m's inequality,
\begin{equation}\label{ineq_2_11}
\frac{x_1^2}{a_1}+\frac{x_2^2}{q_2}+...+\frac{x_n^2}{a_n}\geq\frac{(x_1+x_2+...+x_n)^2}{a_1+a_2+...+a_n}
\end{equation}
for every $x_k\in\RR$ and $a_k>0, \ k\in\{1,2,...,n\}$.
If we make substitutions $a_i=\ds\frac{x_i}{\sqrt{a_i}}$ and
$b_i=\sqrt{a_i}$, for all $i=\{1,2,...,n\}$, in inequality (\ref{ineq_2_10}) we
find a new refinement of Bergstr\"{o}m's inequality, which is given as
follows
\begin{eqnarray*}
&& \left(\sum_{i=1}^n a_i\right)^{-1}\left(\frac{A^2}{M^2}+\frac{2A}{M}\sum_{i=1}^n |x_i|\right) \\
&& \leq \frac{x_1^2}{a_1}+\frac{x_2^2}{a_2}+...+\frac{x_n^2}{a_n}-\frac{(|x_1|+|x_2|+...+|x_n|)^2}{a_1+a_2+...+a_n}\\
&& \leq \left(\sum_{i=1}^n a_i\right)^{-1}\left(\frac{A^2}{m^2}+\frac{2A}{m}\sum_{i=1}^n |x_i|\right),
\end{eqnarray*}
where
$$A=\frac{1}{8}\sqrt{\ds\left(\sum_{i=1}^n \frac{x_i^2}{a_i}\right)\left(\sum_{i=1}^n a_i\right)}\cdot \sum_{i=1}^n \left(\frac{x_i^2}{a_i\ds\sum_{i=1}^n \frac{x_i^2}{a_i}}-\frac{a_i}{\ds\sum_{i=1}^n a_i}\right)^2,$$
$$m=\min_{i=\ov{1,n}}\left\{\frac{x_i^2}{a_i\ds\sum_{i=1}^n \frac{x_i^2}{a_i}},\frac{a_i}{\ds\sum_{i=1}^n a_i}\right\} \ and  \ M=\max_{i=\ov{1,n}}\left\{\frac{x_i^2}{a_i\ds\sum_{i=1}^n \frac{x_i^2}{a_i}},\frac{a_i}{\ds\sum_{i=1}^n a_i}\right\}.$$
\end{itemize}
\end{rema}

\section{APPLICATION TO ARITHMETIC FUNCTIONS}
In the theory of the arithmetic functions \cite{Apo,Nat,SC}, for positive integer $n$,
several important functions have been studied. Among
these we found $\sigma_k(n), \ \tau(n), \ \sigma_k^*(n)$ and
$\tau^*(n)$, where $\sigma_k(n)$ is the sum of $k$th powers of the
divisors of $n$, $\tau(n)$ is the number of divisors of $n$,
$\sigma_k^*(n)$ is the sum of $k$th powers of the unitary divisors
of $n$ and $\tau^*(n)$ is the number of unitary divisors of $n$,
where $k\geq0$.

\begin{teo}
For $n\geq1$ and $k\geq0$, there are the following inequalities
\begin{eqnarray}
\ds\frac{1}{2n\tau(n)}\left[\sigma_{2k}(n)-\left(\frac{\sigma_k(n)}{\tau(n)}\right)^2\right]&\leq&
\frac{\sigma_k(n)}{\tau(n)}-\sqrt{n^k} \nonumber \\
&\leq& \ds\frac{1}{2\tau(n)}\left[\sigma_{2k}(n)-\left(\frac{\sigma_k(n)}{\tau(n)}\right)^2\right] \label{ineq_3_1} 
\end{eqnarray}
and
\begin{eqnarray}
\ds\frac{1}{2n\tau^*(n)}\left[\sigma_{2k}^*
(n)-\left(\frac{\sigma^*_k(n)}{\tau^*(n)}\right)^2\right] &\leq&
\frac{\sigma_k^* (n)}{\tau^*(n)}-\sqrt{n^k} \nonumber \\
&\leq& \ds\frac{1}{2\tau^*(n)}\left[\sigma_{2k}^*(n)-\left(\frac{\sigma_k^*(n)}{\tau^*(n)}\right)^2\right]. \label{ineq_3_2} 
\end{eqnarray}
\end{teo}

\noindent{\it Proof.} If $d_1,d_2,...,d_s$ are the divisors of $n$,
then we take $\alpha_i=\ds\frac{1}{s}$ and $x_i=d_i^k$ in inequality
(\ref{ineq_1_1}).\\
Therefore, we have $m=1, \ M=n$ and $s=\tau(n)$, so inequality (\ref{ineq_1_1})
becomes:
$$\frac{1}{2ns}\sum_{i=1}^s \left(d_i^k -\frac{\sigma_k(n)}{\tau(n)}\right)^2 \leq \frac{\sigma_k(n)}{\tau(n)}-\sqrt{n^k}
\leq\frac{1}{2s}\sum_{i=1}^s\left(d_i^k -\frac{\sigma_k(n)}{\tau(n)}\right)^2.$$
Making simple calculations and taking into account that
$\left(\ds\prod_{i=1}^s d_i^k\right)^{1/s}=\ds\left(\prod_{i=1}^s
d_i\right)^{\frac{k}{s}}=(n^{\frac{s}{2}})^{\frac{k}{s}}=n^{\frac{k}{2}},$
we observe that this inequality is equivalent to inequality (\ref{ineq_3_1}).
Similarly prove that inequality (\ref{ineq_3_2}) is true.

\hfill \qed 

\begin{rema}
Inequality (\ref{ineq_3_2}) improves the inequality
\begin{equation}
\frac{\sigma_k(n)}{\tau(n)}\geq\sqrt{n^k},
\end{equation}
which is due to S. Sivaramakrishnan and C. S. Venkataraman \cite{SC} and
inequality (3.2) improves the inequality
$$\frac{\sigma_k^*(n)}{\tau^*(n)}\geq \sqrt{n^k},$$
which is due to J. S\'{a}ndor and L. T\'{o}th \cite{ST,SC}.
\end{rema}

\section{APPLICATIONS TO OPERATORS}
In this section, we consider bounded linear operators acting on a complex
Hilbert space $\mathcal{H} $. If a bounded linear operator $A$
satisfies $A=A^*$, then $A$ is called a self-adjoint operator. If a
self-adjoint operator $A$ satisfies $\langle x \vert A \vert x
\rangle \geq 0$ for any $\vert x \rangle \in \mathcal{H}$, then $A$
is called a positive operator and denoted by $A\geq 0$.
In addition, $A \geq B$ means $A -B \geq 0$.
We also define operator mean by $A\sharp_{\lambda}B\equiv
A^{1/2}\left(A^{-1/2}BA^{-1/2}\right)^{\lambda}A^{1/2}$ for $\lambda
\in [0,1]$, two invertible positive operator $A$ and $B$ \cite{KA}.
Note that we have the relation $B\sharp_{1-\lambda}A =
A\sharp_{\lambda}B$.
\begin{teo} \label{op_the1}
For $\lambda\in [0,1]$,  two invertible positive operator $A$ and $B$, we have the following relations.
\begin{itemize}
\item[(i)] If $A \leq B$, then we have
\begin{eqnarray}
\frac{\lambda(1-\lambda)}{2} \left(AB^{-1}A-2A+B \right) &\leq& (1-\lambda) A +\lambda B -A \sharp_{\lambda} B \nonumber \\
&\leq & \frac{\lambda(1-\lambda)}{2} \left(BA^{-1}B-2B+A \right). \label{op_the_ineq01}
\end{eqnarray}
\item[(ii)] If $B \leq A$, then we have
\begin{eqnarray}
\frac{\lambda(1-\lambda)}{2}  \left(BA^{-1}B-2B+A \right) &\leq& (1-\lambda) A +\lambda B -A \sharp_{\lambda} B \nonumber \\
&\leq & \frac{\lambda(1-\lambda)}{2}\left(AB^{-1}A-2A+B \right)  .\label{op_the_ineq02}
\end{eqnarray}
\end{itemize}
\end{teo}

{\it Proof:}
We prove (i).
Exchanging $\lambda$ and $1-\lambda$ in the inequalities (1.2), we have
$$
\frac{\lambda(1-\lambda)}{2b}(a-b)^2\leq (1-\lambda)a+\lambda b-a^{1-\lambda}b^{\lambda} \leq \frac{\lambda(1-\lambda)}{2a}(a-b)^2
$$
in the case of $a \leq b$.
Thus we have the inequalities for $0<t \leq 1$:
$$
\frac{\lambda(1-\lambda)}{2}(t-1)^2\leq (1-\lambda)t+\lambda -t^{1-\lambda} \leq \frac{\lambda(1-\lambda)}{2}\left(\sqrt{t}-\frac{1}{\sqrt{t}}\right)^2
$$
putting $t \equiv \frac{a}{b}$.
Thus we have for $0<T\leq I$,
$$
\frac{\lambda(1-\lambda)}{2}(T-1)^2\leq (1-\lambda)T+\lambda -T^{1-\lambda} \leq \frac{\lambda(1-\lambda)}{2}(T^{1/2}-T^{-1/2})^2
$$
by standard operational calculus.
Putting $T=B^{-1/2}AB^{-1/2}$ and then multiplying $B^{1/2}$ from the both sides,
we obtain the desired results.
(ii) can be proven by the similar way to the proof of (i).
\hfill \qed

\begin{rema}
We have $AB^{-1}A-2A+B=A^{1/2}\left(  A^{1/2}B^{-1}A^{1/2}   +\left( A^{1/2}B^{-1}A^{1/2}   \right)^{-1}-2I\right)A^{1/2} \geq 0$,
becasue we have $u+u^{-1}-2=\frac{(u-1)^2}{u}\geq 0$ for scalar $u\geq 0$.
By the similar way, we have  $BA^{-1}B-2B+A \geq 0$.
Thus under the condition $A\leq B$ or $B\leq A$,
the inequalities in Theorem \ref{op_the1} improve the second inequality of the following inequalities (See \cite{FY,Furuichi} for example):
\begin{equation} \label{op_rem_ineq01}
\left\{ (1-\lambda)A^{-1}+\lambda B^{-1} \right\}^{-1} \leq A\sharp_{\lambda} B \leq (1-\lambda) A+\lambda B
\end{equation}
\end{rema}

\begin{cor}
For $\lambda\in (0,1)$,  two invertible positive operator $A$ and $B$, we have the following relations.
\begin{itemize}
\item[(i)] If $A \leq B$, then we have
\begin{eqnarray}
&&A \sharp_{\lambda} B -\left(A \sharp_{\lambda} B \right) \left\{\frac{2}{\lambda(1-\lambda)}\left(B^{-1}AB^{-1}-2B^{-1}+A^{-1}\right)^{-1}+A \sharp_{\lambda} B \right\}\left( A \sharp_{\lambda} B \right) \nonumber\\
  &&\leq \left\{ (1-\lambda) A^{-1} +\lambda B^{-1}\right\}^{-1} \nonumber \\
&&\leq  A \sharp_{\lambda} B -\left(A \sharp_{\lambda} B \right) \left\{\frac{2}{\lambda(1-\lambda)}\left(A^{-1}BA^{-1}-2A^{-1}+B^{-1}\right)^{-1}+A \sharp_{\lambda} B \right\}\left( A \sharp_{\lambda} B \right). \label{op_the_ineq03}
\end{eqnarray}
\item[(ii)] If $B \leq A$, then we have
\begin{eqnarray}
&&A \sharp_{\lambda} B -\left(A \sharp_{\lambda} B \right) \left\{\frac{2}{\lambda(1-\lambda)}\left(A^{-1}BA^{-1}-2A^{-1}+B^{-1}\right)^{-1}+A \sharp_{\lambda} B \right\}\left( A \sharp_{\lambda} B \right)  \nonumber\\
  &&\leq \left\{ (1-\lambda) A^{-1} +\lambda B^{-1}\right\}^{-1} \nonumber \\
&&\leq  A \sharp_{\lambda} B -\left(A \sharp_{\lambda} B \right) \left\{\frac{2}{\lambda(1-\lambda)}\left(B^{-1}AB^{-1}-2B^{-1}+A^{-1}\right)^{-1}+A \sharp_{\lambda} B \right\}\left( A \sharp_{\lambda} B \right). \label{op_the_ineq04}
\end{eqnarray}
\end{itemize}
\end{cor}
{\it Proof:}
Replacing $A$ and $B$ by $A^{-1}$ and $B^{-1}$ in the inequalities (\ref{op_the_ineq01}), respectively and taking the inverse of bothe sides, then
we have the inequalities (\ref{op_the_ineq03}), using $\left(A^{-1} \sharp_{\lambda}B^{-1}\right)^{-1}=  A\sharp_{\lambda}B$ and
$$
\left(  X^{-1}+Y^{-1}  \right)^{-1} = X^{-1} -X^{-1}\left( X^{-1}+Y^{-1}  \right)^{-1}X^{-1}
$$
for invertible positive operators $X$ and $Y$.
The inequalities (\ref{op_the_ineq04}) can be proven by the similar way to the inequalities (\ref{op_the_ineq03}).

\hfill \qed

\begin{rema}
Since $\left(A^{-1}BA^{-1}-2A^{-1}+B^{-1}\right)^{-1}\geq 0$, $\left(B^{-1}AB^{-1}-2B^{-1}+A^{-1}\right)^{-1}\geq 0$ and  $A\sharp_{\lambda}B \geq 0 $,
then the right hand side of the inequalities (\ref{op_the_ineq03}) and (\ref{op_the_ineq04}) are further bounbed from the above by $A\sharp_{\lambda}B $.
Therefore two inequalities (\ref{op_the_ineq03}) and (\ref{op_the_ineq04}) improve the first inequality of the inequalities (\ref{op_rem_ineq01})
under the condition $A \leq B$ or $B \leq A$.
\end{rema}

\begin{cor}
If $0<A\leq B$, then we have
\begin{equation} \label{op_cor_eq01}
3(A-B) + BA^{-1}B-AB^{-1}A \geq 0.
\end{equation}
\end{cor}

The inequality (\ref{op_cor_eq01}) corresponds to the following relation:
$$
0< a \leq b \Rightarrow (b-a)^3 \geq 0
$$
in the commutative case. The inequality (\ref{op_cor_eq01}) can be directly proven by applying the standard operational calculus to
the scalar inequality $(t-1)^3\leq 0$ for $0<t\leq 1$.

\section*{Acknowledgement}
The first author (N.M.) was supported in part by the Romanian
Ministry of Education, Research and Innovation through the PNII Idei
project 842/2008.
 The second author (S.F.) was supported in
part by the Japanese Ministry of Education, Science, Sports and
Culture, Grant-in-Aid for Encouragement of Young Scientists (B),
20740067.


\begin{thebibliography}{10}
\bibitem{Apo}T. M. APOSTOL, {\it Introduction to Analytic Number Theory}, Springer-Verlag, New York, 1976.

\bibitem{CF} D. I. CARTWRIGHT, M. J. FIELD, A refinement of the arithmetic mean-geometric mean inequality, {\it Proc. Amer. Math. Soc.}, Vol.71(1978), pp.36-38.

\bibitem{Furuichi}  S. FURUICHI, On refined Young inequalities and reverse inequalities, {\it J. Math. Ineq.}, Vol.5(2011), pp.21-31.

\bibitem{FM} S. FURUICHI, N. MINCULETE, {\it Alternative reverse inequalities for Young's inequality}, arXiv:1103.1937.

\bibitem{FY} T.FURUTA, M.YANAGIDA, Generalized means and convexity of inversion for positive operators, {\it Amer.Math.Monthly}, Vol.105 (1998),pp.258-259.

\bibitem{HLP} G. H. HARDY, J. E. LITLLWOOD, G. POLYA, {\it Inequalities}, Cambrige University Press, 1934.

\bibitem{KA} F. KUBO, T.ANDO, Means of positive operators, {\it Math. Ann.},Vol.264(1980),pp.205-224.

\bibitem{Nat} M. NATHANSON, {\it Elementary Methods in Number Theory}, Springer, New York, 2006.

\bibitem{NP} C. P. NICULESCU, L.-E. PERSSON, {\it Convex Functions and Their Applications}, CMS Books in Mathematics, Vol. 23, Springer-Verlag, New York, 2006.

\bibitem{Pop} O. T. POP, About Bergstr\"{o}m's inequality, {\it Journal of Mathematical Inequalities}, Vol. 3(2009), pp.213-216.

\bibitem{SC} J. S\'{A}NDOR, B. CRISTICI, {\it Handbook of Number Theory II}, Kluwer Academic Publishers, Dordrecht/Boston/London, 2004.

\bibitem{ST} J. S\'{A}NDOR, L. T\'{O}TH, On certain number-theoretic inequalities, {\it Fib. Quart.} Vol.28 (1990), pp.255-258.

\end{thebibliography}
\end{document}